\pdfoutput=1

\documentclass[12pt]{article}
\usepackage[margin=1in]{geometry}
\usepackage{amsmath}
\usepackage{amsfonts}
\usepackage{amssymb}
\usepackage{amsthm}
\usepackage{esint}
\usepackage{mathtools}
\usepackage{array,makecell,longtable}
\usepackage{titling}
\usepackage{authblk}
\usepackage{commath}
\usepackage{xcolor}
\usepackage{graphicx}
\usepackage{graphbox}
\usepackage[toc,title,page]{appendix}
\usepackage{subcaption}
\usepackage{bbm}
\usepackage{hyperref}
\usepackage{color,soul}
\usepackage{bm}
\usepackage{mathrsfs}

\usepackage{enumitem}
\newtheorem{remark}{Remark}[section]
\newtheorem{theorem}[remark]{Theorem}

\newtheorem{lemma}[remark]{Lemma}

\newtheorem{definition}[remark]{Definition}

\numberwithin{equation}{section}

\title{Determining state space anomalies in mean field games}
\author{Hongyu Liu\thanks{Department of Mathematics, City University of Hong Kong, Hong Kong SAR, China\\ Email address: hongyu.liuip@gmail.com, hongyliu@cityu.edu.hk} \, and Catharine W.K. Lo\thanks{Liu Bie Ju Centre for Mathematical Sciences and Department of Mathematics, City University of Hong Kong, Hong Kong SAR, China\\ Email address: wingkclo@cityu.edu.hk} \vspace{-0.5cm}}

\date{}

\begin{document}

\maketitle

\begin{abstract}
    In this paper, we are concerned with the inverse problem of determining anomalies in the state space associated with the stationary mean field game (MFG) system. We establish novel unique identifiability results for the intrinsic structure of these anomalies in mean field games systems, including their topological structure and parameter configurations, in several general scenarios of practical interest, including traffic flow, market economics and epidemics. To the best of our knowledge, this is the first work that considers anomalies in the state space for the nonlinear coupled MFG system.

    \medskip
		
		

\end{abstract}

\section{Introduction}
\subsection{Problem Setup and Background}
Mean Field Games (MFGs) study the averaged-out behaviours of a large population of agents as the number of agents grows to infinity. The theory of MFGs was pioneered independently by Caines-Huang-Malham\'e \cite{huang2006large} and Lasry-Lions \cite{LasryLions1, LasryLions2} in 2006, and it has since received significant attention and increasing studies in the literature. One of the main characteristics of an MFG is the existence of an adversarial regime, and the Nash equilibrium exists and is unique in this so-called monotone regime.
In its typical formulation, the mean field equilibrium can be characterised by the following MFG system:
\begin{equation}\label{MFG}
    \begin{cases}
		-\partial_t u(x,t) - D \Delta u(x,t) + H(x,\nabla u(x,t))-F(x,t,m(x,t)) = 0,&  \text{ in } \Omega'\times [0,T],\\
		\partial_t m(x,t) - D \Delta m(x,t) - \text{div} (m(x,t) \nabla_p H(x,\nabla u(x,t)) = 0, & \text{ in } \Omega'\times[0,T],\\
        \partial_\nu u(x,t) = \partial_\nu m(x,t)=0 & \text{ on } \partial\Omega'\\
		u(x,T), \quad m(x,0) \quad \text{ given } & \text{ in } \Omega'
    \end{cases}
\end{equation}
in the bounded Lipschitz domain $\Omega'\subset\mathbb{R}^n$ for the Euclidean space $\mathbb{R}^n$ with $n\in\mathbb{N}$. 
Here, $\Delta$ and $\text{div}$ denote the Laplacian and divergent operators with respect to the $x$-variable, respectively, $D(x)>0$ is a coefficient describing the diffusion process, while $H(x,\nabla u) = H(x,p): \Omega'\times\mathbb{R}^n\to\mathbb{R}$ is a nonlinear Hamiltonian. 

In the physical setup, $u$ is the value function of each player; $m$ signifies the population distribution and satisfies the probability density constraint $\int_{\Omega'}m=1$; $F$ is the running cost function which signifies the interaction between the agents and the population; $m_0$ represents the initial population distribution and $G$ signifies the terminal cost. All the functions involved are real valued.

Let $\Omega$ be a closed proper subdomain of $\Omega'$ with a smooth boundary $\partial\Omega$. In the long optimisation time limit $T\to\infty$, the MFG system has an infinite time horizon, and \eqref{MFG} adopts the stationary form
\begin{equation}\label{MFGStat}
    \begin{cases}
		- D \Delta u(x) + H(x,\nabla u(x))+\lambda-F(x,m(x)) = 0,&  \text{ in } \Omega,\\
		- D \Delta m(x) - \text{div} (m(x) \nabla_p H(x,\nabla u(x)) = 0, & \text{ in } \Omega,\\
        u(x) = \psi(x), \quad m(x)=\varphi(x) & \text{ on } \partial\Omega,
    \end{cases}
\end{equation}
where $\lambda$ is a constant that can be determined through the normalisation of $m$. We will consider this system such that there is some topological structure in the system, given by an anomalous inhomogeneity $\omega\Subset\Omega$, which is a bounded Lipschitz domain such that $\Omega\backslash\bar{\omega}$ is connected. Here, by an inhomogeneity, we mean that $H$, $\lambda$ and $F$ has a discontinuity across the boundary of $\omega$. In particular, $H$ has a jump of the form 
\[H(x,\nabla u(x)) = \begin{cases}
    H_1(x,\nabla u(x)) & \text{ if }x\in\omega,\\
    H_0(x,\nabla u(x)) & \text{ otherwise },
\end{cases}
\] with $H_1, H_0$ $C^\gamma$ H\"older continuous for some $\gamma\in(0,1)$ in their respective domains, such that
\[H_1(x,\nabla u(x))\neq H_0(x,\nabla u(x)) \quad \text{ for }x\in\partial\omega.\] A similar form holds for $\lambda$ and $F$. See Definitions \ref{kappaAdmis}--\ref{FAdmis1} and \ref{FDef} in Section \ref{sect:MainResults} for more precise descriptions of these functions. 
This paper aims to quantitatively analyse the geometry of $\omega$ and its intrinsic properties described by $H$ and $F$.

In order to characterise the anomaly $\omega$, we introduce the inverse boundary problem given by the following measurement map of a single pair of Cauchy data:
\begin{equation}\label{MeasureMap}
    \mathcal{M}_{\omega,H,F}  := \left((\psi,\varphi),(\partial_\nu u, \partial_\nu m)_{\partial\Omega}\right),\, \psi,\varphi \text{ fixed} \to \omega,H,F,
\end{equation}
where $\nu$ is the exterior unit normal vector to $\partial\Omega$. In the physical context, $\omega$ signifies the support of the anomalous polyhedral inclusion in the state space $\Omega$. Hence, the inverse problem \eqref{MeasureMap} is concerned with recovering the location and shape of this anomaly, as well as its parameter configuration. It is also referred to as the inverse inclusion problem in the theory of inverse problems.

For the inverse inclusion problem \eqref{MeasureMap}, we mainly consider its unique identifiability issue. That is, we aim at establishing the sufficient conditions under which $\omega$ can be uniquely determined by the measurement map $\mathcal{M}_{\omega,H,F}$ in the sense that if two admissible inclusions $(\omega_j, H_j, F_j)$, $j = 1, 2$, produce the same boundary measurement, i.e. $\mathcal{M}_{\omega_1,H_1,F_1}(u,m)=\mathcal{M}_{\omega_2,H_2,F_2}(u,m)$ associated with a fixed $(u,m)$, then one has $(\omega_1, H_1, F_1)=(\omega_2, H_2, F_2)$. 

Formally, our main result can be summarised as follows: 
\begin{theorem}\label{thm:formal}
    Consider a certain general scenario where there exists an anomalous polyhedral inclusion $\omega\Subset\Omega$, such that $H$, $\lambda$ and $F$ has a jump discontinuity given by their adherence to certain admissibilility properties. Then $\omega$ and its parameter configurations $H$ and $F$ are uniquely determined by the boundary measurement $\mathcal{M}_{\omega,H,F}$.
\end{theorem}
The details of the anomalous inclusion $\omega$ and the admissible class of the unknowns $H$, $\lambda$ and $F$ will be provided in Section \ref{sect:MainResults}. The main result is contained in Theorems \ref{thm:Statkappa} and \ref{thm:StatF}. 

We shall establish sufficient conditions such that Theorem \ref{thm:formal} holds in a certain general setup. In particular, we shall provide general characterisations of $H$, $F$ and $\omega$, and connect them to specific applications. The major novelty that distinguishes our inverse problem study from most of the existing ones lies at the following three aspects. First, we consider inverse inclusion problems for coupled nonlinear parabolic partial differential equations (PDEs) with discontinuous parameters, which have never been considered elsewhere. To address this problem, we innovate the method of microlocal analysis for corner singularities that make full use of the intrinsic structure of the MFG system. Furthermore, this method allows us to overcome the technical constraints of positivity and probability density, rendering our results more physically realistic. Our study opens up a new direction of research on inverse problems for mean field games with many potential developments and physical applications. More discussion on these aspects shall be given in Section \ref{sect:discuss}.

\subsection{Technical Developments and Discussion }\label{sect:discuss}

Developed in recent years, mean field game (MFG) is a powerful framework for modelling large-scale interactive systems involving numerous agents. These rational agents aim to optimise their individual objectives while taking into account the collective behavior of the entire population. These nonlinear coupled PDE systems were first studied in the seminal works of Lasry-Lions \cite{LasryLions1,LasryLions2,lasry2007mean}, and independently by Caines-Huang-Malham\'e \cite{huang2006large}. Since then, MFGs have witnessed significant growth and have become a prominent research area (see, for instance, \cite{cardaliaguet2010notes}, \cite{cardaliaguet2015weak}, \cite{CardaliaguetGraber2015mean}, \cite{CardaliaguetGraberPorrettaTonon2015second}, \cite{CarmonaDelarue2018_1}, \cite{CarmonaDelarue2018_2}, \cite{cardaliaguet2019master}, \cite{CardaliaguetPorretta2020Lectures}, \cite{Ambrose2022existence}, \cite{meszaros2023mean} and the references therein), with the stationary case being studied in \cite{LasryLions1}, \cite{Cirant2016MFGStationary}, \cite{ABC2017StationaryMFG}, \cite{FerreiraGomesTada2019StationaryMFGDirichlet}, and \cite{Osborne2022StationaryMFGDirichlet} among many others. These foundational works paved the way for a deeper exploration of MFGs in various contexts, including economics, engineering, and social sciences.

While the well-posedness of the MFG system is well understood, the inverse problems for MFGs are far less studied in the literature. Related works only started appearing in the previous year, beginning with numerical results in \cite{chow2022numerical}, \cite{ding2022mean}, \cite{mou2022numerical} and \cite{klibanov2023convexification}, and later, theoretical results in \cite{LiuMouZhang2022InversePbMeanFieldGames}, \cite{LiuZhang2022-InversePbMFG}, \cite{klibanov2023lipschitz}, \cite{klibanov2023mean1}, \cite{klibanov2023mean2}, \cite{LiuZhangMFG3}, \cite{klibanov2023holder}, \cite{liu2023stability}, \cite{imanuvilov2023lipschitz1}, \cite{imanuvilov2023unique}, \cite{klibanov2023coefficient1}, 
\cite{klibanov2023coefficient2}, \cite{imanuvilov2023global} and \cite{imanuvilov2023lipschitz2}. However, all these works dealt with the unique identifiability results for functions in the MFG system, such as the solution, the running cost, the total cost or the Hamiltonian.

In this work, we instead focus on anomalies in the state space, such as domains with singular corners. We aim to recover such anomalies from the behaviour of the solution at the boundary of the domain. This anomaly is portrayed as a discontinuity in the Hamiltonian $H$ and the running cost $F$, which had always been assumed to be continuous in the previous works on MFGs mentioned above. Such an inverse problem for MFG has a wide range of applications. For instance, when MFGs are applied in the study of traffic flow or other forms of congestion management \cite{MFGCar2}, \cite{MFGCrowd}, \cite{MFGCar}, \cite{MFGCrowd+Econs} we can make use of boundary measurements of the traffic to detect an obstruction \cite{MFGStateConstraintsTrafficObstruct}, \cite{ghattassi2023nonseparable}. An example is where we can measure the flow of cars in and out of a city, to detect any bridge collapse or car accident \cite{MFGCarAccident}. Such applications can also be extended to other forms of flow optimisation in engineering, including short circuits in smart grids \cite{MFGSmartGrid2} and breaks in supply chains \cite{MFGSupplyChain}. Moreover, MFGs can also be applied in the development of autonomous driving \cite{MFGAutoCar1,MFGAutoCar2,MFGAutoCar3}, and our work can help in the detection of corners and walls or other types of anomalous inhomogeneities. On the other hand, when MFGs are applied in the context of economics to model market dynamics, price formation or resource allocation \cite{MFGCrowd+Econs}, \cite{CarmonaDelarue2018_1}, our work can allow for the recovery of external influences limiting the price, including governmental influences, limit ups/downs or irreversible investments \cite{MFGIrreversibleInvestJump}. In the context of social sciences, MFGs can also model epidemic spreading \cite{MFGDisease2} and disease propagation \cite{MFGDisease1}, as well as social network information spread \cite{MFGSocialNetwork}. In these cases, there might be gaps in the state space due to deserted areas of the country, vaccinated parts of the population, or internet censorship. Our work can be applied to detect such anomalies. 
Such problems of detecting anomalies have previously been considered in the context of forward problems by some authors (see, for instance, \cite{MFG_B+C+C_RegimeSwitchJump}, \cite{dianetti2023ergodic} and \cite{MFGSingularity} and the references therein), but to the best of our knowledge, there has not yet been any results in the realm of inverse problems.

Since we are working with the MFG coupled system, we take inspiration from the general use of higher order linearisation techniques for solving nonlinear equations in \cite{LinLiuLiuZhang2021-InversePbSemilinearParabolic-CGOSolnsSuccessiveLinearisation}, \cite{LiuMouZhang2022InversePbMeanFieldGames}, and adopt a similar approach to solve the inverse problem in a nonlinear coupled system. This enables us to break up the nonlinear system of equations into several uncoupled linear or nonlinear equations. 
However, it should be noted that in the works \cite{LinLiuLiuZhang2021-InversePbSemilinearParabolic-CGOSolnsSuccessiveLinearisation}, \cite{LiuMouZhang2022InversePbMeanFieldGames}, or other similar works, the authors considered the dense subspace of the solution space, given by all possible oscillating complex geometric solutions around the trivial
solution $(u,m) \equiv (0, 0)$. As a result, the solutions are not always non-negative, which is important for MFGs where $m$ represents the non-negative probability density. Consequently, even if the results in \cite{LiuMouZhang2022InversePbMeanFieldGames} are mathematically rigorous, this work is unrealistic from a physical perspective. Indeed, this technical constraint significantly increases the difficulty of the inverse problem, since one would need to construct suitable ``probing modes" to fulfil this constraint. Several works have tried to address this technical constraint, such as in \cite{LiuZhangMFG3}, and also in the stationary case in \cite{ding2023determining} (see also \cite{liu2023determining} where the authors also treated the issue of positivity, but in the case of biological systems). 

In this paper, we make use of the intrinsic features of our inverse inclusion problem to develop a novel approach to ensure the positivity constraint on $m$, while effectively tackling the MFG inverse problems. 
A crucial element of this approach is that we are not required to consider all the possible solutions. On the other hand, we can fix the boundary measurements we need. Therefore, we can control these boundary values, so that the solution $m$ satisfies the positivity constraint, as discussed in the first paragraphs of the proofs of Theorems \ref{thm:Statkappa} and \ref{thm:StatF} in Sections \ref{sect:ThmStatKappaPf} and \ref{sect:ThmStatFPf} respectively. This novel approach can effectively tackle the constraints mentioned above, and may have potential for further uses in situations when only a finite number of measurements are required, such as in stability studies.

Yet, it should be noted that since we are making use of the method of higher order linearisation, an infinite number of measurements is required. This is different from the minimal number of measurements required in \cite{DiaoFeiLiuWang-2022-Semilinear-Shape-Corners}, and is a technical constraint arising from nonlinearity.

Moreover, we note that in this case, unlike the works of \cite{LinLiuLiuZhang2021-InversePbSemilinearParabolic-CGOSolnsSuccessiveLinearisation} and \cite{LiuMouZhang2022InversePbMeanFieldGames} and similar works, we do not require $F$ to be regular or analytic. In fact, we consider the case where $F$ has a jump across the boundary of $\omega$. This is the main difference which separates our work from previous works on nonlinear elliptic equations.

In summary, the major contributions of this work are as follows: 
Firstly, we analysed corner singularities at the microlocal level in dimensions $n\geq2$. Such analysis has never been done for nonlinear coupled systems, and has only been previously achieved for $n\leq3$ in a different context.
Secondly, we established uniqueness results for the MFG system in recovering anomalies/obstacles. These results are highly interesting, in particular in the following three aspects. 
First, to our best knowledge, this is the first result in the literature concerning the shape determination of anomalous inclusions in the MFG system, which is a nonlinear coupled system that has not been previously considered elsewhere. This is highly relevant to physical applications, as seen in previous works we have cited above. 
Second, in our MFG inverse problem, we consider the case where the running cost $F$ is not continuous and has a jump across the boundary of $\omega$. This is another novelty.
Third, our results are achieved by specific choice of measurements, which helps to overcome the technical constraints of the problem as we discussed previously. 
We are confident that the mathematical techniques described in this work can be used to address more inverse problems connected to coupled nonlinear PDEs in various contexts.

The rest of the paper is organized as follows. In Section \ref{sect:Prelims}, we give the basic notations and present the well-posedness of the forward MFG problem in Section \ref{sect:Notations}. This is followed by a description of the geometric structures of the anomalies we are considering, in Section \ref{sect:geometry}. With these initial setups in place, we formulate our main results mathematically in Section \ref{sect:MainResults}. We then develop the higher order linearisation technique in Section \ref{sect:Linearise}, which is necessary for the proofs of the main results in Section \ref{sect:MainProof}. The main results corresponding to the Hamiltonian $H$ and the running cost $F$ are then provided in Sections \ref{sect:ThmStatKappaPf} and \ref{sect:ThmStatFPf} respectively.
	
\section{Preliminaries}\label{sect:Prelims}
	
\subsection{Notations and Basic Setting}\label{sect:Notations}

Let $\mathcal{P}$ stand for the set of Borel probability measures on $\mathbb{R}^n$, and $\mathcal{P}(\Omega')$ stand for the set of Borel probability measures on $\Omega'$. Let $m\in\mathcal{P}(\Omega')$ denote the population distribution of the agents and $u(x, t):\Omega'\times [0, T]\to \mathbb{R}$ denote the value function of each player. We restrict our study to MFGs with quadratic Hamiltonians given by the following form:

\begin{equation}\label{MFGQuadraticStat}
    \begin{cases}
        -D \Delta u(x) + \frac{1}{2}\kappa(x) |\nabla u(x)|^2 + \lambda - F(x,m(x)) = 0 &\quad \text{in }\Omega',\\
        -D \Delta m(x) - \text{div}(\kappa(x) m(x)\nabla u(x)) = 0  &\quad \text{in }\Omega',\\
        \partial_\nu u(x)=\partial_\nu m(x)=0 &\quad \text{on }\partial\Omega', \\
        u(x) = \psi(x), \quad m(x)=\varphi(x) & \quad \text{on } \partial\Omega,
    \end{cases}
\end{equation} 
for a Lipschitz domain $\Omega'\subset\mathbb{R}^n$, and $F:\Omega'\times\mathcal{P}(\Omega')\to\mathbb{R}$. Then, it is well-known (see, for instance, Theorem 2.1 of \cite{LasryLions1}) that
\begin{theorem}\label{ForwardPbMainThm}
    Suppose $F$ is Lipschitz and bounded in $\Omega'$. Then there exists $\lambda\in\mathbb{R}$, $u\in C^2(\Omega')$ and $m\in W^{1,p}(\Omega')$ for any $p<\infty$ solving \eqref{MFGQuadraticStat}.
\end{theorem}

For $k\in\mathbb{N}$ and $0<\alpha<1$, recall that the H\"older space $C^{k+\alpha}(\overline{\Omega})$ is defined as the subspace of $C^{k}(\overline{\Omega})$ such that $\phi\in C^{k+\alpha}(\overline{\Omega})$ if and only if $D^l\phi$ exist and are H\"older continuous with exponent $\alpha$ for all $l=(l,l,\dots,l_n)\in\mathbb{N}^n$ with $|l|\leq k$, where $D^l:=\partial^{l}_{x}\partial^{l}_{x}\cdots\partial^{l_n}_{x_n}$ for $x=(x,\cdots,x_n)$. The corresponding norm is defined by
\[\norm{\phi}_{C^{k+\alpha}(\overline{\Omega})}:=\sum_{|l|\leq k}\norm{D^l\phi}_\infty+\sum_{|l|= k}\sup_{x\neq y}\frac{|D^l\phi(x)-D^l\phi(y)|}{|x-y|^\alpha}\] for the $L^\infty$ sup norm $\norm{\cdot}_\infty$.

\subsection{Geometrical Setup}\label{sect:geometry}

Consider the convex conic cone $\mathcal{S}\subset\Omega$ with apex $x_c$ and axis $v_c$ and opening angle $2\theta_c\in(0,\pi)$, defined by
\[\mathcal{S}:=\{y\in\Omega:0\leq\angle (y-x_c,v_c)\leq\theta_c,\theta_c\in(0,\pi/2)\}.\] Define the truncated conic cone by
\[\mathcal{S}_h:=\mathcal{S}\cap B_h,\] where $B_h:=B_h(x_c)$ is an open ball contained in $\Omega$ centred at $x_c$ with radius $h>0$. Observe that both $\mathcal{S}$ and $\mathcal{S}_h$ are Lipschitz domains. 

We first consider some asymptotics for a CGO solution we will be using. 
\begin{lemma}\label{wCGOlem}
Let $w$ be the solution to 
\begin{equation}\label{CGOEq}
- \Delta w(x) = 0\text{ in }\Omega',
\end{equation}
of the form 
\begin{equation}\label{CGO}
w = e^{\tau (\xi + i\xi^\perp)\cdot (x-x_c)}
\end{equation} 
such that $\xi\cdot\xi^\perp=0$, $\xi,\xi^\perp\in\mathbb{S}^{n-1}$. 
Then there exists a positive number $\rho$ depending on $\mathcal{S}_h$ satisfying 
\begin{equation}\label{CGOCond}
-1 < \xi\cdot\widehat{(x-x_c)}\leq -\rho < 0\quad\text{ for all } x\in \mathcal{S}_h,
\end{equation}
where $\hat{x} = \frac{x}{|x|}$. Moreover, for sufficiently large $\tau$, \cite{DiaoFeiLiuWang-2022-Semilinear-Shape-Corners}
\begin{equation}\label{CGOEst1}
\left|\int_{\mathcal{S}_h}w\right|\geq C_{\mathcal{S}_h}\tau^{-n} +\mathcal{O}\left(\frac{1}{\tau}e^{-\frac{1}{2}\rho h \tau}\right),
\end{equation}
\begin{equation}\label{CGOEst2}
\left|\int_{\mathcal{S}_h}|x-x_c|^\alpha w\right|\lesssim \tau^{-(\alpha+n)}+\frac{1}{\tau}e^{-\frac{1}{2}\rho h \tau} \quad\forall \alpha>0,
\end{equation}
\begin{equation}\label{CGOEst3}
\norm{w}_{H^1(\partial\mathcal{S}_h)}\lesssim (2\tau^2+1)^{\frac12}e^{-\rho h\tau},
\end{equation}
\begin{equation}\label{CGOEst4}
\norm{\partial_\nu w}_{L^2(\partial\mathcal{S}_h)}\lesssim \tau e^{-\rho h \tau}.
\end{equation}
Here, we use the symbol `` $\lesssim$" to denote that the inequality holds up to a constant which is independent of $\tau$.
\end{lemma}

\begin{proof}
    We first remark that since $\xi\perp\xi^\perp$, clearly $w$ satisfies $D \Delta w=0$ in $\mathbb{R}^n$, and in particular, in $\Omega'$.

    Furthermore, \eqref{CGOCond} is easily satisfied with an appropriate choice of $\xi$, since $\mathcal{S}_h$.

    Next, we recall from Lemma 2.2 of \cite{DiaoFeiLiuYang-2022-EM-Shape-Corners} that, for some fixed $\alpha>0$ and $0<\delta<e$,
    \begin{equation}\label{LaplaceTransformIdentity}
        \int_0^\delta r^\alpha e^{-\mu r}\,dr = \frac{\Gamma(\alpha+1)}{\mu^{\alpha+1}} + \int_\delta^\infty r^\alpha e^{-\mu r}\,dr,
    \end{equation} 
    where $\mu\in\mathbb{C}$ and $\Gamma$ is the Gamma function. Moreover, if the real part of $\mu$, denoted by $\mathscr{R}\mu$, is such that $\mathscr{R}\mu \geq \frac{2\alpha}{e}$, then $r^\alpha\leq e^{\mathscr{R}\mu r/2}$, and hence 
    \begin{equation}\label{LaplaceTransformIneq}
        \left|\int_\delta^\infty r^\alpha e^{-\mu r}\,dr\right| \leq \frac{2}{\mathscr{R}\mu}e^{-\mathscr{R}\mu\delta/2}.
    \end{equation}
    With this, we write $x-x_c=(x_1,\dots,x_n)\in\mathbb{R}^n$ in polar coordinates:
    \[x-x_c= \left( \begin{array}{l}
    r\cos\theta_1 \\
    r\sin\theta_1\cos\theta_2 \\
    r\sin\theta_1\sin\theta_2\cos\theta_3 \\
    \vdots \\
    r\sin\theta_1\cdots\sin\theta_{n-2}\cos\varphi \\
    r\sin\theta_1\cdots\sin\theta_{n-2}\sin\varphi
    \end{array} \right)^T,\]
    where $\theta_i\in[0,\pi]$ for $i=1,\dots,n-2$ and $\varphi\in[0,2\pi)$, with Jacobian
    \[r^{n-1}\sin^{n-2}(\theta_1)\sin^{n-3}(\theta_2)\cdots\sin(\theta_{n-2})\,dr\,d\theta_1 \,d\theta_2 \cdots\,d\theta_{n-2}\,d\varphi.\]

    Since $\mathcal{S}_h$ has an opening angle of $2\theta_c$, by \eqref{LaplaceTransformIdentity}, we have 
    \begin{multline}\label{CGOEstProofEq}
    \int_{\mathcal{S}_h}e^{\tau (\xi + i\xi^\perp)\cdot (x-x_c)} \\ = I_1 + \int_0^{2\pi} \left(\int_0^{\theta_c} \cdots \int_0^{\theta_c} I_2 \sin^{n-2}(\theta_1)\sin^{n-3}(\theta_2)\cdots\sin(\theta_{n-2})\,d\theta_1\,d\theta_2\cdots\,d\theta_{n-2} \right) \,d\varphi, 
    \end{multline}
    where 
    \begin{equation}\label{CGOEstProofEqI1}I_1 := \int_0^{2\pi} \int_0^{\theta_c} \cdots \int_0^{\theta_c}\frac{\Gamma(n)}{\tau^n \left((\xi + i\xi^\perp) \cdot \widehat{(x-x_c)}\right)^n} \sin^{n-2}(\theta_1)\cdots\sin(\theta_{n-2})\,d\theta_1\cdots\,d\theta_{n-2}\,d\varphi\end{equation} 
    and 
    \begin{equation}\label{CGOEstProofEqI2}I_2 := \int_h^\infty r^{n-1} e^{\tau r (\xi + i\xi^\perp)\cdot \widehat{(x-x_c)}} \,dr.\end{equation}
    By the integral mean value theorem, we have that 
    \begin{align}\label{CGOEstProofEq1}
        I_1 & = \frac{\Gamma(n)}{\tau^n} \int_0^{2\pi} \frac{1}{\left((\xi + i\xi^\perp) \cdot \widehat{(x-x_c)}(\varphi,\theta_\zeta)\right)^n} \,d\varphi\int_0^{\theta_c} \sin^{n-2}(\theta_1) \,d\theta_1 \cdots \int_0^{\theta_c}\sin(\theta_{n-2})\,d\theta_{n-2} \nonumber \\
        & = \frac{2\pi\Gamma(n)C_{\theta_c}}{\tau^n} \frac{1}{\left((\xi + i\xi^\perp) \cdot \widehat{(x-x_c)}(\varphi_\zeta,\theta_\zeta)\right)^n} 
    \end{align}
    for some constant $C_{\theta_c}$ such that $0<C_{\theta_c}<1$. Since $\xi \cdot \widehat{(x-x_c)} > -1$ by \eqref{CGOCond}, this gives that 
    \begin{equation}\label{CGOEstProofEq2}
    |I_1| \geq \frac{2\pi\Gamma(n)C_{\theta_c}}{\tau^n} \frac{1}{2^{n/2}}. 
    \end{equation}

    On the other hand, for sufficiently large $\tau$, by \eqref{LaplaceTransformIneq},
    \begin{equation}|I_2| = \left| \int_h^\infty r^{n-1} e^{r \tau (\xi + i\xi^\perp)\cdot \widehat{(x-x_c)}} \,dr \right| \leq \frac{2}{\rho\tau} e^{-\frac{1}{2}\rho h \tau}.\end{equation} 
    Consequently, we have that 
    \begin{multline}
    \left| \int_0^{2\pi} \left(\int_0^{\theta_c} \cdots \int_0^{\theta_c} I_2 \sin^{n-2}(\theta_1)\sin^{n-3}(\theta_2)\cdots\sin(\theta_{n-2})\,d\theta_1\,d\theta_2\cdots\,d\theta_{n-2} \right) \,d\varphi \right| 
    \\
    \leq \int_0^{2\pi} \left(\int_0^{\theta_c} \cdots \int_0^{\theta_c} |I_2| \,d\theta_1\,d\theta_2\cdots\,d\theta_{n-2} \right) \,d\varphi  \leq \frac{4 \pi \theta_c^{n-2}}{\rho \tau } e^{-\frac{1}{2}\rho h \tau}.
    \end{multline}
    Combining this with \eqref{CGOEstProofEq2}, we obtain \eqref{CGOEst1} with $C_{\mathcal{S}_h} = \frac{\pi\Gamma(n)C_{\theta_c}}{2^{n/2-1}}$.

    Next, for \eqref{CGOEst2}, utilising \eqref{LaplaceTransformIdentity}, we obtain higher exponents for $r$ in \eqref{CGOEstProofEq}, i.e. 
    \begin{multline*}
    \int_{\mathcal{S}_h}(x-x_c)^\alpha e^{\tau (\xi + i\xi^\perp)\cdot (x-x_c)} \\ = I_{3} + \int_0^{2\pi} \left(\int_0^{\theta_c} \cdots \int_0^{\theta_c} I_{4} \sin^{n-2}(\theta_1)\sin^{n-3}(\theta_2)\cdots\sin(\theta_{n-2})\,d\theta_1\,d\theta_2\cdots\,d\theta_{n-2} \right) \,d\varphi, 
    \end{multline*}
    where 
    \begin{align*}& |I_{3}| \\
    & := \left|\int_0^{2\pi} \int_0^{\theta_c} \cdots \int_0^{\theta_c}\frac{\Gamma(n+\alpha)}{\tau^{n+\alpha} \left((\xi + i\xi^\perp) \cdot \widehat{(x-x_c)}\right)^{n+\alpha}} \sin^{n-2}(\theta_1)\cdots\sin(\theta_{n-2})\,d\theta_1 \cdots\,d\theta_{n-2}\,d\varphi \right|\\
    & = \left| \frac{2\pi\Gamma(n+\alpha)C_{\theta_c}}{\tau^{n+\alpha}} \frac{1}{\left((\xi + i\xi^\perp) \cdot \widehat{(x-x_c)}(\varphi_\zeta,\theta_\zeta)\right)^{n+\alpha}} \right| \\
    & \lesssim\frac{1}{\tau^{n+\alpha}} 
    \end{align*} 
    as in the procedure to obtain \eqref{CGOEstProofEq1}--\eqref{CGOEstProofEq2}, 
    and 
    \begin{equation*}|I_{4}| := \left| \int_h^\infty r^{n+\alpha-1} e^{r \tau (\xi + i\xi^\perp)\cdot \widehat{(x-x_c)}} \,dr \right| \leq \frac{2}{\rho\tau} e^{-\frac{1}{2}\rho h \tau},\end{equation*} 
    which directly leads to \eqref{CGOEst2}.

    At the same time, 
    \begin{equation}\norm{w}_{H^1(\partial\mathcal{S}_h)} = \left( \norm{w}_{L^2(\partial\mathcal{S}_h)}^2 + \norm{\tau (\xi + i\xi^\perp)w}_{L^2(\partial\mathcal{S}_h)}^2 \right)^{\frac{1}{2}} \leq (2\tau^2+1)^{\frac{1}{2}} \norm{w}_{L^2(\partial\mathcal{S}_h)} \end{equation} 
    by the Cauchy-Schwarz inequality, since $\xi\in\mathbb{S}^{n-1}$. 
    But, applying the polar coordinate transformation on $w$, we have that 
    \begin{equation}\label{CGOEstProofEq3}\norm{w}_{L^2(\partial\mathcal{S}_h)} = \left( \int_0^{2\pi} \int_0^{\theta_c} \cdots \int_0^{\theta_c} e^{2\tau h \xi \cdot \widehat{(x-x_c)}} \,d\theta_1 \cdots \, d\theta_{n-2} \, d\varphi \right)^{\frac{1}{2}} \leq (2\pi \theta_c^{n-2})^{\frac{1}{2}} e^{-\rho h\tau}\end{equation} 
    by \eqref{CGOCond}. This gives \eqref{CGOEst3}.

    Finally, 
    \begin{equation}\norm{\partial_\nu w}_{L^2(\partial\mathcal{S}_h)} \leq \norm{\nabla w}_{L^2(\partial\mathcal{S}_h)} = \norm{\tau (\xi + i\xi^\perp)w}_{L^2(\partial\mathcal{S}_h)} \leq \sqrt{2}\tau \norm{w}_{L^2(\partial\mathcal{S}_h)} \end{equation} once again by the Cauchy-Schwarz inequality. Using \eqref{CGOEstProofEq3} then gives the desired result \eqref{CGOEst4}.
\end{proof}

This generalises Lemma 2.1 of \cite{DiaoFeiLiuWang-2022-Semilinear-Shape-Corners} to any dimension $n\geq2$, and for general corners with apex at $x_c$ not necessarily the origin.

\subsection{Main Results}\label{sect:MainResults}

With these definitions in place, we can now state the main results of this paper. Let $\Omega$ be a closed proper subdomain of $\Omega'$ with smooth boundary $\partial\Omega$. Consider a bounded Lipschitz polyhedral domain $\omega\Subset\Omega$ such that $\Omega\backslash\bar{\omega}$ is connected and each corner of $\omega$ is a convex cone $\mathcal{S}$, as defined in Subsection \ref{sect:geometry}. We define the following admissibility conditions for $\kappa$, $\lambda$ and $F$.
\begin{definition}\label{kappaAdmis}
    We say that $\kappa(x):\Omega'\to\mathbb{R}$ is admissible, denoted by $\kappa\in\mathcal{B}$, if $\kappa(x)$ is of the form 
    \[\kappa(x) = \kappa_0(x) + (\kappa_1(x)-\kappa_0(x))\chi_\omega, \quad x\in\Omega'\]
    such that $\kappa_1$ and $\kappa_0$ are $C^\gamma$ H\"older-continuous for some $\gamma\in(0,1)$ with respect to $x\in \omega$ and $x\in \Omega\backslash \omega$ respectively, and in an open neighbourhood $U$ of $\partial \omega$, $\kappa_1(U\cap \omega) \neq \kappa_0(U\cap(\Omega\backslash \omega))$.
\end{definition}

\begin{definition}\label{lambdaAdmis}
    We say that the constant $\lambda$ is admissible, denoted by $\lambda\in\mathcal{C}$, if $\lambda$ is of the form 
    \[\lambda = \lambda_0 + (\lambda_1-\lambda_0)\chi_\omega\]
    such that $\lambda_0$ and $\lambda_1$ are (possibly different) constants.
\end{definition}

Next, we introduce the admissible class of locally analytic functions.
\begin{definition}\label{FAdmis1}
    Let $E$ be a compact subset of $\Omega\subset\mathbb{R}^n$. We say that $U(x,z):E\times\mathbb{C}\to\mathbb{C}$ is admissible, denoted by $U\in\mathcal{A}_E$, if:
    \begin{enumerate}[label=(\alph*)]
        \item The map $z\mapsto U(\cdot,z)$ is holomorphic with value in $C^{2+\alpha}(E)$ for some $\alpha\in(0,1)$, 
        \item $U(x,\cdot)$ is $C^\sigma$-continuous with respect to $x\in E$ for some $\sigma\in(0,1)$. 
        \item $U(x,0)=\lambda_U$ for all $x\in \mathbb{R}^n$.
    \end{enumerate}
    It is clear that if $U$ satisfies these two conditions, $U$ can be expanded into a power series
    \[U(x,z)=\lambda_U+\sum_{m=1}^\infty U^{m}(x)\frac{z^m}{m!},\]
    where $U^{m}(x)=\frac{\partial^m }{\partial z^m}U(x,0)\in C^{2+\alpha}(E)$.
\end{definition}

This admissibility condition is imposed a priori on $U$, by extending these functions of real variables to the complex plane with respect to the $z$-variable, given by $\tilde{U}$, and assume that they are holomorphic as functions of the complex variables $z$. Then, $U$ is simply the restriction of $\tilde{U}$ to the real line. Furthermore, since the image of $U$ is in $\mathbb{R}$, we can assume that the series expansion of $\tilde{U}$ is real-valued.

We remark that this is different from the admissibility condition considered in previous works on MFGs, including \cite{LiuMouZhang2022InversePbMeanFieldGames}, \cite{LiuZhangMFG3} and \cite{ding2023determining}, where the authors assume that $F$ is analytic in the entire real space $\mathbb{R}^n$. Our case is more general, and allows for singularities at the boundary $\partial E$ of compact subsets $E$ of $\Omega$.

Observe that if $F\in\mathcal{A}$ with the appropriate $\lambda$, then $(u,m)=(0,0)$ is a solution to \eqref{MFGQuadraticStat}. In fact, it is easy to check that $\left(u,\exp\left(\int \kappa(x)\nabla u(x)\,dx\right)\right)$ is also a solution to \eqref{MFGQuadraticStat}, for some $u$.  
With this admissibility condition on $F$, we are able to improve the result of Theorem \ref{ForwardPbMainThm} and obtain a stronger local regularity of the solution $(u,m)$ to the MFG system \eqref{MFGQuadraticStat}, using the implicit functions theorem for Banach spaces. We omit the proof here, since it is simply a minor modification of Theorem 2.1 of \cite{ding2023determining}.

\begin{theorem}\label{ForwardResultStat}
Suppose that $F\in\mathcal{A}_E$ for the compact subset $E$ of $\Omega$. Then, there exist constants $\delta>0$ and $C>0$ such that for any 
    \[\psi(x),\varphi(x)\in B_{\delta}(C^{2+\alpha}(\partial\Omega)) :=\{h\in C^{2+\alpha}(\partial\Omega): \|\psi(x)\|_{C^{2+\alpha}(\Omega)}\leq\delta \}, \]
    the MFG system \eqref{MFGQuadraticStat} has a solution $(u,m)\in [C^{2+\alpha}(E)]^2$ which satisfies
    \begin{equation}
		\|(u,m)\|_{ C^{2+\alpha}(\overline{E})}:= \|u\|_{C^{2+\alpha}(\overline{E})}+ \|m\|_{C^{2+\alpha}(\overline{E})}\leq C(\|\psi(x)\|_{ C^{2+\alpha}(\partial\Omega)}+\|\varphi(x)\|_{ C^{2+\alpha}(\partial\Omega)}).
    \end{equation}
    Furthermore, the solution $(u,m)$ is unique within the class
    \begin{equation}
		\{ (u,m)\in  C^{2+\alpha}(\overline{E})\times C^{2+\alpha}(\overline{E}): \|(u,m)\|_{ C^{2+\alpha}(\overline{E})}\leq C\delta \}.
    \end{equation}		
\end{theorem}

With this local well-posedness in hand, we can define the admissibility class of the running cost $F$.

\begin{definition}\label{FDef}
    We say that $F(x,m):\Omega'\times\mathcal{P}(\Omega')\to\mathbb{R}$ is admissible, denoted by $F\in\mathcal{D}$, if 
    $F$ is of the form 
    \[F(x) = F_0(x) + (F_1(x)-F_0(x))\chi_\omega, \quad x\in\Omega',\]
    such that $F_1\in\mathcal{A}_\omega$ with $F_1(x,0)=\lambda_1$, $F_0\in\mathcal{A}_{\Omega\backslash\omega}$ with $F_0(x,0)=\lambda_0$, and for some $\ell\in\mathbb{N}$, the $\ell$-th Taylor coefficient of $F_1$ and $F_0$, denoted by $F^{(\ell)}_1(x)$ and $F^{(\ell)}_0(x)$ respectively, are $C^\sigma$ H\"older-continuous for some $\sigma\in(0,1)$ with respect to $x\in \omega$ and $x\in \Omega\backslash \omega$ respectively, and in an open neighbourhood $U$ of $\partial \omega$, $F^{(\ell)}_1(U\cap \omega) \neq F^{(\ell)}_0(U\cap(\Omega\backslash \omega))$.
\end{definition}

With the above admissibility conditions, then, microlocally in a corner $\mathcal{S}_h$ of $\omega$, we have that $(u_1,m_1)$ and $(u_0,m_0)$ satisfies \begin{equation}\label{MFGQuadraticStat1}
    \begin{cases}
        -D \Delta u_1(x) + \frac{1}{2}\kappa_1(x) |\nabla u_1(x)|^2 + \lambda_1 - F_1(x,m_1(x)) = 0 &\quad \text{in }\mathcal{S}_h,\\
        -D \Delta m_1(x) - \text{div}(\kappa_1(x) m_1(x)\nabla u_1(x)) = 0  &\quad \text{in }\mathcal{S}_h,\\
        \partial_\nu u_1(x)=\partial_\nu m_1(x)=0 &\quad \text{on }\partial\Omega',\\
        u_1(x) = \psi(x), \quad m_1(x)=\varphi(x) & \quad\text{on } \partial\Omega.
    \end{cases}
\end{equation} and 
\begin{equation}\label{MFGQuadraticStat0}
    \begin{cases}
        -D \Delta u_0(x) + \frac{1}{2}\kappa_0(x) |\nabla u_0(x)|^2 + \lambda_0 - F_0(x,m_0(x)) = 0 &\quad \text{in }B_h\backslash\mathcal{S}_h,\\
        -D \Delta m_0(x) - \text{div}(\kappa_0(x) m_0(x)\nabla u_0(x)) = 0  &\quad \text{in }B_h\backslash\mathcal{S}_h,\\
        \partial_\nu u_0(x)=\partial_\nu m_0(x)=0 &\quad \text{on }\partial\Omega',\\
        u_0(x) = \psi(x), \quad m_0(x)=\varphi(x) & \quad\text{on } \partial\Omega
    \end{cases}
\end{equation}
respectively, such that 
\begin{equation}
    \begin{cases}
        \partial_\nu u_0(x)=\partial_\nu u_1(x), \quad \partial_\nu m_0(x) = \partial_\nu m_1(x) &\quad \text{on }\partial\mathcal{S}_h\backslash\partial B_h,\\
        u_0(x) = u_1(x), \quad m_0(x)=m_1(x) & \quad\text{on } \partial\mathcal{S}_h\backslash\partial B_h.
    \end{cases}
\end{equation}

We next state our main results.

\begin{theorem}\label{thm:Statkappa}
    Suppose that $\kappa\in\mathcal{B}$, $\lambda\in\mathcal{C}$ and $F\in\mathcal{A}_\Omega$. Then $\omega$ and $\kappa$ are uniquely determined by the boundary measurement $\mathcal{M}_{\omega,H,F}$.
\end{theorem}

\begin{theorem}\label{thm:StatF}
    Suppose that $\lambda\in\mathcal{C}$, $F\in\mathcal{D}$ . Then $\omega$ and $F$ are uniquely determined by the boundary measurement $\mathcal{M}_{\omega,H,F}$. 
\end{theorem}

\section{Higher-Order Linearisation}\label{sect:Linearise}

We next develop a high-order linearisation scheme of the MFG system \eqref{MFGQuadraticStat} with respect to $\psi(x)$ and $\varphi(x)$, in the case that $\tilde{F}\in\mathcal{A}_E$ with $\tilde{F}(x,0)-\lambda=0$.  This is possible due to the infinite differentiability of the MFG system \eqref{MFGQuadraticStat} in $E\Subset\Omega$ with respect to $\psi,\varphi$, as a consequence of Theorem \ref{ForwardResultStat}. 
Let 
\[\psi(x;\varepsilon)=\sum_{\ell=1}^N\varepsilon^\ell f_\ell \quad \text{ and }\quad \varphi(x;\varepsilon)=\sum_{\ell=1}^N\varepsilon^\ell g_\ell\quad \text{ for }x\in \partial\Omega, \]
where 
\[f_\ell\in \{f \in C^{2+\alpha}(\partial\Omega):\norm{f}_{C^{2+\alpha}(\partial\Omega)} < \delta\}\] and \[g_\ell\in \{g \in C^{2+\alpha}(\partial\Omega):g>0\text{ and }\norm{g}_{C^{2+\alpha}(\partial\Omega)} < \delta\}\] for a small enough $\delta$ so that Theorem \ref{ForwardResultStat} holds
and 
$\varepsilon\in\mathbb{R}_+$ with 
$|\varepsilon|$ small enough. The choice of $\varphi$ to be positive can be easily ensured, since we are choosing our boundary measurements.

By Theorem \ref{ForwardResultStat}, there exists a unique solution $(u(x;\varepsilon),m(x;\varepsilon))$ of \eqref{MFGQuadraticStat}. Since $\tilde{F}(x,0)-\lambda=0$, when $\varepsilon=0$, we have $(u(x;0),m(x;0))= (0,0)$. 

Let $S$ be the solution operator of \eqref{MFGQuadraticStat} defined in Theorem \ref{ForwardResultStat}. Then there exists a bounded linear operator $L$ from $\mathcal{H}:=[B_{\delta}( C^{2+\alpha}(\partial\Omega))]^2$ to $[C^{2+\alpha}(E)]^2$ such that
\begin{equation}
	\lim\limits_{\|(\psi,\varphi)\|_{\mathcal{H}}\to0}\frac{\|S(\psi,\varphi)-S(u(x;0),m(x;0))- L(\psi,\varphi)\|_{[C^{2+\alpha}(E)]^2}}{\|(\psi,\varphi)\|_{\mathcal{H}}}=0,
\end{equation} 
where $\|(\psi,\varphi)\|_{\mathcal{H}}:=\|\psi\|_{ C^{2+\alpha}(\partial\Omega)      }+\|\varphi\|_{C^{2+\alpha}(\partial\Omega)}$.

Now we consider $\varepsilon=0$.
Then it is easy to check that $L(\psi,\varphi)|_{\varepsilon=0}$ is the solution map of the following system, as in \cite{LiuMouZhang2022InversePbMeanFieldGames}, \cite{LiuZhang2022-InversePbMFG} ,\cite{LiuZhangMFG3} or \cite{ding2023determining}: 
\begin{equation}\label{Linear1}
    \begin{cases}
		-D \Delta u^{(1)}(x)= F^{(1)}(x)m^{(1)}(x)  & \text{ in }  \Omega',\\
		-D \Delta m^{(1)}(x)=0 & \text{ in }  \Omega',\\
        \partial_\nu u^{(1)}(x)=\partial_\nu m^{(1)}(x)=0 & \text{ on }\partial\Omega',\\
		u^{(1)}(x)=f_1(x), \quad m^{(1)}(x)=g_1(x) & \text{ on } \partial\Omega.
    \end{cases}
\end{equation}
This system is called the first-order linearization system.

In the following, we define
\begin{equation}\label{eq:ld1}
 (u^{(1)}, m^{(1)} ):=L(\psi,\varphi)\Big|_{\varepsilon=0}. 
 \end{equation}
For notational convenience, we write
\begin{equation}\label{eq:ld2}
u^{(1)}=\partial_{\varepsilon}u(x;\varepsilon)|_{\varepsilon=0}\quad\text{and}\quad m^{(1)}=\partial_{\varepsilon}m(x;\varepsilon)|_{\varepsilon=0}.
\end{equation}
In our subsequent discussion, we will employ these notations to streamline the presentation, and their intended significance will be evident within the given context.

Next, for the second-order linearisation, we consider 
\[u^{(2)}:=\partial_{\varepsilon}^2u|_{\varepsilon=0},\quad 	m^{(2)}:=\partial_{\varepsilon}^2m|_{\varepsilon=0}.\]
Then, we have the second-order linearisation as follows:
\begin{equation}\label{Linear2}
    \begin{cases}
		-D \Delta u^{(2)}(x)+\kappa(x)|\nabla u^{(1)}|^2= F^{(1)}m^{(2)}(x)+F^{(2)}[m^{(1)}(x)]^2& \text{ in } \Omega'\\
		-D \Delta m^{(2)}(x)= 2\text{div} (m^{(1)}\kappa(x)\nabla u^{(1)}) &\text{ in } \Omega',\\
        \partial_\nu u(x)=\partial_\nu m(x)=0 & \text{ on }\partial\Omega',\\
		u^{(2)}(x)=2f_2(x), \quad m^{(2)}(x)=2g_1(x) & \text{ on } \partial\Omega.\\
    \end{cases}  	
\end{equation}

Notice that the non-linear terms of the system \eqref{Linear2} depend on the first-order linearised system \eqref{Linear1}. Therefore, to consider higher order Taylor coefficients of $F$, we consider, for $\ell\in\mathbb{N}$, 
\[u^{(\ell)}=\partial_{\varepsilon}^\ell u|_{\varepsilon=0}, \quad	m^{(\ell)}=\partial_{\varepsilon}^\ell m|_{\varepsilon=0}.\]
In this manner, we produce a series of elliptic systems that will be used once more when there is a discontinuity in the higher order Taylor coefficient $F^{(\ell)}$.

\section{Proof of Main Results}\label{sect:MainProof}

We next proceed to prove the main Theorems \ref{thm:Statkappa} and \ref{thm:StatF}. Before that, we begin first with a main auxiliary theorem. 

\begin{theorem}\label{AuxThm}
    Let $\omega\Subset\Omega$ be the bounded Lipschitz domain such that $\Omega\backslash\bar{\omega}$ is connected, with conic corner $\mathcal{S}_h$. For $P\in L^2(\Omega)$, suppose that $q_1,q_0$ are $C^\gamma$ H\"older-continuous for some $\gamma\in(0,1)$ with respect to $x\in \omega$, such that $q_1\neq q_0$ in $\mathcal{S}_h$. For any function $h\in C^1(\Omega)$, consider the following system of equations for $v_j\in C^{2+\alpha}(\mathcal{S}_h)$, $j=1,0$:
    \begin{equation}\label{AuxThmEq}
    \begin{cases}
		-D \Delta v_j(x)+ P + q_j(x)h(x)= 0& \quad \text{ in } \mathcal{S}_h,\\
		\partial_\nu v_1(x)=\partial_\nu v_0(x) &\quad \text{ on }\partial\mathcal{S}_h\backslash\partial B_h,\\
        v_1(x)=v_0(x) &\quad \text{ on }\partial\mathcal{S}_h\backslash\partial B_h.
    \end{cases}  	
\end{equation}
Then, it must hold that \[h(x_c) =0,\] where $x_c$ is the apex of $\mathcal{S}_h$.
\end{theorem}

\begin{proof}

Taking the difference between the $j$-th equation, we have, denoting $\tilde{v}=v_1-v_0$
\begin{equation}\label{Linear2jUDiff}
    \begin{cases}
		-D \Delta \tilde{v}(x)+(q_1(x)-q_0(x))h(x)= 0& \quad \text{ in } \mathcal{S}_h\\
		\partial_\nu \tilde{v}(x)=0 & \quad \text{ on }\partial\mathcal{S}_h\backslash\partial B_h,\\
		\tilde{v}(x)=0 & \quad \text{ on } \partial\mathcal{S}_h\backslash\partial B_h.
    \end{cases}  	
\end{equation}
Multiplying this by the solution $w$ to \eqref{CGOEq} given by \eqref{CGO} and integrating in the truncated cone $\mathcal{S}_h$, we have, by Green's formula, 
\begin{equation}\label{k1DiffEq}
    \int_{\mathcal{S}_h}(q_1(x) -q_0(x))h(x)  w(x) \,dx = - \int_{\partial \mathcal{S}_h} w(x) \partial_\nu \tilde{v}(x) \, d\sigma + \int_{\partial \mathcal{S}_h} \tilde{v}(x) \partial_\nu w(x) \, d\sigma.
\end{equation}

Since $q_j(x)\in C^\gamma(\mathcal{S}_h)$, $\gamma\in(0,1)$ and $h \in C^1(\Omega)$ for $j=1,0$, the function 
\begin{equation}\label{barkappa}\bar{q}:=(q_1(x) -q_0(x))h(x)\in C^\gamma(\mathcal{S}_h).\end{equation} 
Therefore, we can expand $\bar{q}$ as follows:
\[\bar{q}=\bar{q}(x_c)+\delta\bar{q},\quad |\delta\bar{q}|\leq \norm{\bar{q}}_{C^\gamma(\mathcal{S}_h)}|x-x_c|^\gamma,\] 
where \[\bar{q}(x_c)=(q_1(x_c) -q_0(x_c))h(x_c).\] Thus, the left-hand-side of \eqref{k1DiffEq} can be expanded as
\begin{align*}
    \int_{\mathcal{S}_h}\bar{q}  w(x) \,dx = \bar{q}(x_c)\int_{\mathcal{S}_h} w(x) \,dx + \int_{\mathcal{S}_h}\delta\bar{q}w \,dx,
\end{align*}
where
\begin{align*}
    \left|\int_{\mathcal{S}_h}\delta\bar{q}w \,dx\right| \leq \norm{\bar{q}}_{C^\gamma(\mathcal{S}_h)} \int_{\mathcal{S}_h}|x-x_c|^\gamma |w| \,dx    
\end{align*}
for 
\begin{equation}\label{PfEq5}\norm{\bar{q}}_{C^\gamma(\mathcal{S}_h)}\leq \norm{h}_{C^\gamma(\mathcal{S}_h)}\sup_{\mathcal{S}_h}\left|q_1 - q_0\right| + \norm{q_1 -q_0}_{C^\gamma(\mathcal{S}_h)}\sup_{\mathcal{S}_h} h\leq C\end{equation} for some constant $C$ depending on the norm of $h$. 

On the other hand, the right-hand-side of \eqref{k1DiffEq} can be analysed as follows: By the Cauchy-Schwarz inequality and the trace theorem, 
\begin{align*}
    \left|\int_{\partial \mathcal{S}_h} w(x) \partial_\nu \tilde{v}(x) \, d\sigma\right| & \leq \norm{w}_{H^{\frac12}(\partial \mathcal{S}_h)}\norm{\partial_\nu \tilde{v}}_{H^{-\frac12}(\partial \mathcal{S}_h)} \\
    & \leq C \norm{w}_{H^1(\partial \mathcal{S}_h)}\norm{\tilde{v}}_{H^1(\mathcal{S}_h)}
    \\
    & \lesssim (2\tau^2+1)^{\frac12}e^{-\rho h\tau},
\end{align*}
by \eqref{CGOEst3}, while
\begin{align*}
    \left|\int_{\partial \mathcal{S}_h} \tilde{v}(x) \partial_\nu w(x) \, d\sigma\right| & \leq \norm{\tilde{v}}_{L^2(\partial \mathcal{S}_h)}\norm{\partial_\nu w}_{L^2(\partial \mathcal{S}_h)} \\
    & \leq C \norm{\tilde{v}}_{H^1(\partial \mathcal{S}_h)}\norm{\partial_\nu w}_{L^2(\partial\mathcal{S}_h)}
    \\
    & \lesssim \sqrt{2}\tau e^{-\rho h\tau}.
\end{align*}
by \eqref{CGOEst4}. Note that the (trace or Sobolev) constants $C>0$ here may be different, and different from that in \eqref{PfEq5}.

Combining these estimates and \eqref{CGOEst2}, we have
\begin{align*}
    \bar{q}(x_c)\int_{\mathcal{S}_h} w(x) \,dx & \lesssim \int_{\mathcal{S}_h}|x-x_c|^\gamma |w| \,dx + (2\tau^2+1)^{\frac12}e^{-\zeta h\tau} + \sqrt{2}\tau e^{-\zeta h\tau} \\
    & \lesssim \tau^{-(\gamma+n)}+\frac{1}{\tau}e^{-\frac{1}{2}\rho h \tau} + (2\tau^2+1)^{\frac12}e^{-\rho h\tau} + \sqrt{2}\tau e^{-\rho h\tau}.
\end{align*}
At the same time, applying the estimate \eqref{CGOEst1}, we have 
\begin{align*}
    \bar{q}(x_c) \left[C_{\mathcal{S}_h}\tau^{-n} +\mathcal{O}\left(\frac{1}{\tau}e^{-\frac{1}{2}\rho h \tau}\right)\right] \lesssim \tau^{-(\gamma+n)}+\frac{1}{\tau}e^{-\frac{1}{2}\rho h \tau} + (1+\tau)e^{-\rho h\tau}
\end{align*}
when $\tau\to\infty$.

Multiplying by $\tau^n $ on both sides and letting $\tau\to\infty$, we have that $\bar{q}(x_c)=0$. Since $q_1(x_c) \neq q_0(x_c)$, by the definition of $\bar{q}$ in \eqref{barkappa}, it must be that $h(x_c)=0$.

\end{proof}

\subsection{A Corner Singularity due to Different Hamiltonians}\label{sect:ThmStatKappaPf}

We next proceed to prove Theorem \ref{thm:Statkappa}.

\begin{proof}[Proof of Theorem \ref{thm:Statkappa}]
Since $\psi,\varphi$ is fixed in \eqref{MeasureMap}, we can take $\varphi(x)=\sum_{\ell=1}^N\varepsilon^\ell g_\ell$ such that $g_1=g_2=0$ and $g_3>0$, and $\psi(x)=\sum_{\ell=1}^2\varepsilon^\ell f_\ell$ such that the solution $u^{(1)}$ of \eqref{Linear1} corresponding to $f_1$ has non-zero gradient on $\partial\omega$. This ensures that the results are physically meaningful with $m>0$, as discussed in Section \ref{sect:discuss}.

Now, suppose on the contrary that there exists $\omega_1$ and $\omega_2$, such that $\omega_1\neq\omega_2$. Since $\omega_1,\omega_2$ are polyhedrons, it must be that there exists a corner $\mathcal{S}$ of $\omega_1$ such that $\mathcal{S}\Subset\Omega\backslash\overline{\omega_2}$. We also assume that for $\omega_i$, $\kappa\in\mathcal{B}$ with the functions $\kappa_0$ and $\kappa_i$ for $i=1,2$ respectively. Then, in the neighbourhood $B_h$ of $\mathcal{S}_h$, it holds that 
\begin{equation}\label{MFGQuadraticStatExtCornerkappa}
    \begin{cases}
        -D \Delta u_1(x) + \frac{1}{2}\kappa_1(x) |\nabla u_1(x)|^2 + \lambda -  F(x,m_1(x)) = 0 &\quad \text{in }\mathcal{S}_h,\\
        -D \Delta m_1(x) - \text{div}(\kappa_1(x) m_1(x)\nabla u_1(x)) = 0  &\quad \text{in }\mathcal{S}_h,\\
        -D \Delta u_0(x) + \frac{1}{2}\kappa_0(x) |\nabla u_0(x)|^2 + \lambda - F(x,m_0(x)) = 0 &\quad \text{in }B_h\backslash\mathcal{S}_h,\\
        -D \Delta m_0(x) - \text{div}(\kappa_0(x) m_0(x)\nabla u_0(x)) = 0  &\quad \text{in }B_h\backslash\mathcal{S}_h,\\
        \partial_\nu u_1(x)=\partial_\nu u_0(x), \quad \partial_\nu m_1(x)=\partial_\nu m_0(x) &\quad \text{on }\partial\mathcal{S}_h\backslash\partial B_h,\\
        u_1(x)=u_0(x), \quad m_1(x)=m_0(x) &\quad \text{on }\partial\mathcal{S}_h\backslash\partial B_h.
    \end{cases}
\end{equation}

For $F\in\mathcal{A}_\Omega$ with $F(x,0)=\lambda$, we conduct higher order linearisation separately for $(u_1,m_1)$ and $(u_0,m_0)$ in $\mathcal{S}_h$ and $B_h\backslash\mathcal{S}_h$ respectively. By the choice of $\varphi$ with $g_1=0$, it is straightforward to observe from the first order linearisation \eqref{Linear1} that $m^{(1)}_1=m^{(1)}_0\equiv0$, and subsequently $m^{(2)}_1=m^{(2)}_0\equiv0$ from the second order linearisation \eqref{Linear2} with $g_2=0$. It can also be easily seen that as a result, when $\mathcal{M}_{\omega_1,H_1,F}=\mathcal{M}_{\omega_2,H_2,F}$, $u^{(1)} := u^{(1)}_0 + (u^{(1)}_1 - u^{(1)}_0)\chi_\omega \in C^{2+\alpha}(\Omega)$ is such that $u^{(1)}_0 = u^{(1)}_1$ on $\partial\omega$, since they satisfy 
\begin{equation}
    \begin{cases}
		-D \Delta u^{(1)}_1= 0  & \quad\text{ in }  \mathcal{S}_h,\\
        -D \Delta u^{(1)}_0= 0  & \quad\text{ in }  B_h\backslash\mathcal{S}_h,\\
        \partial_\nu u_1^{(1)}(x)=\partial_\nu u_0^{(1)}(x) &\quad \text{ on }\partial\mathcal{S}_h\backslash\partial B_h,\\
        u_1^{(1)}(x)=u_0^{(1)}(x) &\quad \text{ on }\partial\mathcal{S}_h\backslash\partial B_h,\\
        u_1^{(1)}(x)=u_0^{(1)}(x) = f_1 &\quad \text{ on }\partial\Omega,
    \end{cases}
\end{equation}
from the first order linearisation \eqref{Linear1}. 

On the other hand, $u^{(2)}_j$ satisfies
\begin{equation}\label{Linear2jU}
    \begin{cases}
		-D \Delta u^{(2)}_1(x)+\kappa_1(x)|\nabla u^{(1)}|^2= 0& \quad \text{ in } \mathcal{S}_h,\\
        -D \Delta u^{(2)}_0(x)+\kappa_0(x)|\nabla u^{(1)}|^2= 0& \quad \text{ in } B_h\backslash\mathcal{S}_h,\\
		\partial_\nu u_1^{(2)}(x)=\partial_\nu u_0^{(2)}(x) &\quad \text{ on }\partial\mathcal{S}_h\backslash\partial B_h,\\
        u_1^{(2)}(x)=u_0^{(2)}(x) &\quad \text{ on }\partial\mathcal{S}_h\backslash\partial B_h.
    \end{cases}  	
\end{equation}
Since $u^{(1)}$ is fixed and smooth, this is simply an elliptic equation. Therefore, we can apply the unique continuation principle for elliptic equations \cite{KochTataru2001UCPElliptic} to obtain 
\[
    \begin{cases}
		-D \Delta u^{(2)}_1(x)+\kappa_1(x)|\nabla u^{(1)}|^2 = -D \Delta u^{(2)}_0(x)+\kappa_0(x)|\nabla u^{(1)}|^2 = 0& \quad \text{ in } \mathcal{S}_h,\\
		\partial_\nu u_1^{(2)}(x)=\partial_\nu u_0^{(2)}(x) &\quad \text{ on }\partial\mathcal{S}_h\backslash\partial B_h,\\
        u_1^{(2)}(x)=u_0^{(2)}(x) &\quad \text{ on }\partial\mathcal{S}_h\backslash\partial B_h,
    \end{cases} 
\]
which is simply \eqref{AuxThmEq} with $v_j=u^{(2)}_j$, $P=0$, $q_j=\kappa_j$ and $h=|\nabla u^{(1)}|^2$. Since $\kappa\in\mathcal{B}$ and $u\in C^2(\Omega')$ by Theorem \ref{ForwardPbMainThm}, the assumptions of Theorem \ref{AuxThm} are satisfied, and we have that $h(x_c)=|\nabla u^{(1)}(x_c)|^2=0$, i.e. $\nabla u^{(1)}(x_c)=0$, for the apex $x_c$ of $\mathcal{S}$. 

Yet, we have fixed $\psi(x)$ such that $u^{(1)}$ has non-zero gradient on $\partial\omega_1\supset\partial\mathcal{S}$. Thus, we arrive at a contradiction, and it must be that $\omega_1=\omega_2$ and we have uniquely determined $\omega$.

Moreover, when $\mathcal{M}_{\omega_1,H_1,F}=\mathcal{M}_{\omega_2,H_2,F}$, we have that the solution sets $(u_1,m_1),(u_2,m_2)$ satisfy 
\begin{equation}
    \begin{cases}
        -D \Delta u_1(x) + \frac{1}{2}\kappa_1(x) |\nabla u_1(x)|^2 + \lambda - F(x,m_1(x)) = 0 &\quad \text{in }\mathcal{S}_h,\\
        -D \Delta m_1(x) - \text{div}(\kappa_1(x) m_1(x)\nabla u_1(x)) = 0  &\quad \text{in }\mathcal{S}_h,\\
        -D \Delta u_2(x) + \frac{1}{2}\kappa_2(x) |\nabla u_2(x)|^2 + \lambda - F(x,m_2(x)) = 0 &\quad \text{in }\mathcal{S}_h,\\
        -D \Delta m_2(x) - \text{div}(\kappa_2(x) m_2(x)\nabla u_2(x)) = 0  &\quad \text{in }\mathcal{S}_h,\\
        \partial_\nu u_1(x)=\partial_\nu u_2(x), \quad \partial_\nu m_1(x)=\partial_\nu m_2(x) &\quad \text{on }\partial\mathcal{S}_h\backslash\partial B_h,\\
        u_1(x)=u_2(x), \quad m_1(x)=m_2(x) &\quad \text{on }\partial\mathcal{S}_h\backslash\partial B_h.
    \end{cases}
\end{equation}
Repeating the argument above, we have that $m^{(1)}_1=m^{(1)}_2=m^{(2)}_1=m^{(2)}_2\equiv0$ in $\Omega$, $u^{(1)}_1 = u^{(1)}_2$ in $\mathcal{S}_h$, and $u^{(2)}_j$ satisfies
\begin{equation}\label{PfkappaEq5}
    \begin{cases}
		-D \Delta u^{(2)}_1(x)+\kappa_1(x)|\nabla u^{(1)}|^2 = -D \Delta u^{(2)}_2(x)+\kappa_2(x)|\nabla u^{(1)}|^2 = 0& \quad \text{ in } \mathcal{S}_h,\\
		\partial_\nu u_1^{(2)}(x)=\partial_\nu u_2^{(2)}(x) &\quad \text{ on }\partial\mathcal{S}_h\backslash\partial B_h,\\
        u_1^{(2)}(x)=u_2^{(2)}(x) &\quad \text{ on }\partial\mathcal{S}_h\backslash\partial B_h.
    \end{cases} 
\end{equation}
Then, once again by Theorem \ref{AuxThm}, if $\kappa_1(x_c)\neq\kappa_2(x_c)$, we have that $\nabla u^{(1)}(x_c)=0$. This once again contradicts the choice of $\psi$, so it must be that $\kappa_1(x_c)=\kappa_2(x_c)$ and we have uniquely determined $\kappa$.

Observe that in this second part of the proof for the uniqueness of $\kappa$, it is not necessary to apply the unique continuation principle, since the equations in \eqref{PfkappaEq5} for $u^{(2)}_1$ and $u^{(2)}_2$ are both in $\mathcal{S}_h$.

\end{proof}

This extends Theorem 2.2 of \cite{LiuZhangMFG3} to more general Hamiltonians which are discontinuous across the boundary $\partial\omega$.

Note that in the proof, we only require that the first two Taylor coefficients of $m$ are zero on the boundary of $\Omega$. This allows us to ensure the physical relevance of the problem, by picking positive higher order Taylor coefficients, to ensure that $m$ is positive.

\subsection {A Corner Singularity due to Different Running Costs $F$}\label{sect:ThmStatFPf}

Finally, we prove Theorem \ref{thm:StatF}.

\begin{proof}[Proof of Theorem \ref{thm:StatF}]
Since $\psi,\varphi$ is fixed in \eqref{MeasureMap}, we can take $\psi(x)=0$ and $\varphi(x)=\sum_{\ell=1}^N\varepsilon^\ell g_\ell$, such that the corresponding solution $m^{(1)}$ to $g_1>0$ of \eqref{Linear1} satisfies $m^{(1)}\neq0$ on $\partial\omega$. In particular, the positivity of $g_1$ ensures the positivity of $m$, which makes the problem physically meaningful, as we discussed in Section \ref{sect:discuss}.

As in the proof of Theorem \ref{thm:Statkappa}, we suppose on the contrary that there exists $\omega_1$ and $\omega_2$, such that $\omega_1\neq\omega_2$. Since $\omega_1,\omega_2$ are polyhedrons, it must be that there exists a corner $\mathcal{S}$ of $\omega_1$ such that $\mathcal{S}\Subset\Omega\backslash\overline{\omega_2}$. We also assume that for $\omega_i$, $F\in\mathcal{D}$ with the functions $F_0$ and $F_i$ for $i=1,2$ respectively. Then, in the neighbourhood $B_h$ of $\mathcal{S}_h$, it holds that 
\begin{equation}\label{MFGQuadraticStatExtCornerF}
    \begin{cases}
        -D \Delta u_1(x) + \frac{1}{2}\kappa(x) |\nabla u_1(x)|^2 + \lambda_1 - F_1(x,m_1(x)) = 0 &\quad \text{in }\mathcal{S}_h,\\
        -D \Delta m_1(x) - \text{div}(\kappa(x) m_1(x)\nabla u_1(x)) = 0  &\quad \text{in }\mathcal{S}_h,\\
        -D \Delta u_0(x) + \frac{1}{2}\kappa(x) |\nabla u_0(x)|^2 + \lambda_0 - F_0(x,m_0(x)) = 0 &\quad \text{in }B_h\backslash\mathcal{S}_h,\\
        -D \Delta m_0(x) - \text{div}(\kappa(x) m_0(x)\nabla u_0(x)) = 0  &\quad \text{in }B_h\backslash\mathcal{S}_h,\\
        \partial_\nu u_1(x)=\partial_\nu u_0(x), \quad \partial_\nu m_1(x)=\partial_\nu m_0(x) &\quad \text{on }\partial\mathcal{S}_h\backslash\partial B_h,\\
        u_1(x)=u_0(x), \quad m_1(x)=m_0(x) &\quad \text{on }\partial\mathcal{S}_h\backslash\partial B_h.
    \end{cases}
\end{equation}

For $F_1\in\mathcal{A}_\omega$ and $F_0\in\mathcal{A}_{\Omega\backslash\omega}$ with $F_1(x,0)=\lambda_1$ and $F_0(x,0)=\lambda_0$, we can conduct higher order linearisation separately for $(u_1,m_1)$ and $(u_0,m_0)$ in $\mathcal{S}_h$ and $B_h\backslash\mathcal{S}_h$ respectively. 

We first show the case when $\ell=1$. When $\mathcal{M}_{\omega_1,H,F_1}=\mathcal{M}_{\omega_2,H,F_2}$, we have that $m^{(1)} := m^{(1)}_0 + (m^{(1)}_1 - m^{(1)}_0)\chi_\omega \in C^{2+\alpha}(\Omega)$ with $m^{(1)}_0 = m^{(1)}_1$ on $\partial\omega$, since they satisfy 
\begin{equation}
    \begin{cases}
		-D \Delta m^{(1)}_1= 0  & \quad\text{ in }  \mathcal{S}_h,\\
        -D \Delta m^{(1)}_0 = 0  & \quad\text{ in }  B_h\backslash\mathcal{S}_h,\\
        \partial_\nu m_1^{(1)}(x)=\partial_\nu m_0^{(1)}(x) &\quad \text{ on }\partial\mathcal{S}_h\backslash\partial B_h,\\
        m_1^{(1)}(x)=m_0^{(1)}(x) &\quad \text{ on }\partial\mathcal{S}_h\backslash\partial B_h,\\
        m_1^{(1)}(x)=m_0^{(1)}(x) = g_1 &\quad \text{ on }\partial\Omega,
    \end{cases}
\end{equation}
from the first order linearisation \eqref{Linear1}. 

On the other hand, $u^{(1)}_j$ satisfies
\begin{equation}\label{Linear1m}
    \begin{cases}
		-D \Delta u^{(1)}_1= F^{(1)}_1(x)m^{(1)}(x)  & \quad\text{ in }  \mathcal{S}_h,\\
        -D \Delta u^{(1)}_0= F^{(1)}_0(x)m^{(1)}(x)  & \quad\text{ in }  B_h\backslash\mathcal{S}_h,\\
        \partial_\nu u_1^{(1)}(x)=\partial_\nu u_0^{(1)}(x) &\quad \text{ on }\partial\mathcal{S}_h\backslash\partial B_h,\\
        u_1^{(1)}(x)=u_0^{(1)}(x) &\quad \text{ on }\partial\mathcal{S}_h\backslash\partial B_h.
    \end{cases}
\end{equation}
Once again, observe that this is simply an elliptic equation, so we can apply the unique continuation principle for elliptic equations \cite{KochTataru2001UCPElliptic} to obtain 
\[
    \begin{cases}
		-D \Delta u^{(1)}_1 - F^{(1)}_1(x)m^{(1)}(x) = -D \Delta u^{(1)}_0 - F^{(1)}_0(x)m^{(1)}(x) = 0  & \quad\text{ in }  \mathcal{S}_h,\\
        \partial_\nu u_1^{(1)}(x)=\partial_\nu u_0^{(1)}(x) &\quad \text{ on }\partial\mathcal{S}_h\backslash\partial B_h,\\
        u_1^{(1)}(x)=u_0^{(1)}(x) &\quad \text{ on }\partial\mathcal{S}_h\backslash\partial B_h,
    \end{cases}
\]
which is simply \eqref{AuxThmEq} with $v_j=u^{(1)}_j$, $P=0$, $q_j=F^{(1)}_j$ and $h=m^{(1)}$. Since $F\in\mathcal{D}$ and $m^{(1)}\in C^{2+\alpha}(\mathcal{S}_h)$ by Theorem \ref{ForwardResultStat}, the assumptions of Theorem \ref{AuxThm} are satisfied with $F^{(1)}_1(x_c)\neq F^{(1)}_0(x_c)$, and we have that $h(x_c)=m^{(1)}(x_c)=0$ for the apex $x_c$ of $\mathcal{S}$. 

Yet, we have fixed $\varphi(x)$ such that $m^{(1)}\neq0$ on $\partial\omega_1\supset\partial\mathcal{S}$. Thus, we arrive at a contradiction, and it must be that $\omega_1=\omega_2$ and we have uniquely determined $\omega$.

Moreover, when $\mathcal{M}_{\omega_1,H,F_1}=\mathcal{M}_{\omega_2,H,F_2}$, we have that the solution sets $(u_1,m_1),(u_2,m_2)$ satisfy 
\begin{equation}\label{PfF12}
    \begin{cases}
        -D \Delta u_1(x) + \frac{1}{2}\kappa(x) |\nabla u_1(x)|^2 + \lambda_1 - F_1(x,m_1(x)) = 0 &\quad \text{in }\mathcal{S}_h,\\
        -D \Delta m_1(x) - \text{div}(\kappa(x) m_1(x)\nabla u_1(x)) = 0  &\quad \text{in }\mathcal{S}_h,\\
        -D \Delta u_2(x) + \frac{1}{2}\kappa(x) |\nabla u_2(x)|^2 + \lambda_2 - F_2(x,m_2(x)) = 0 &\quad \text{in }\mathcal{S}_h,\\
        -D \Delta m_2(x) - \text{div}(\kappa(x) m_2(x)\nabla u_2(x)) = 0  &\quad \text{in }\mathcal{S}_h,\\
        \partial_\nu u_1(x)=\partial_\nu u_2(x), \quad \partial_\nu m_1(x)=\partial_\nu m_2(x) &\quad \text{on }\partial\mathcal{S}_h\backslash\partial B_h,\\
        u_1(x)=u_2(x), \quad m_1(x)=m_2(x) &\quad \text{on }\partial\mathcal{S}_h\backslash\partial B_h.
    \end{cases}
\end{equation}
Repeating the argument above, we have that $m^{(1)}_1 = m^{(1)}_2$ in $\mathcal{S}_h$, and $u^{(1)}_j$ satisfies
\begin{equation}
    \begin{cases}
		-D \Delta u^{(1)}_1 - F^{(1)}_1(x)m^{(1)}(x) = -D \Delta u^{(1)}_2 - F^{(1)}_2(x)m^{(1)}(x) = 0  & \quad\text{ in }  \mathcal{S}_h,\\
        \partial_\nu u_1^{(1)}(x)=\partial_\nu u_2^{(1)}(x) &\quad \text{ on }\partial\mathcal{S}_h\backslash\partial B_h,\\
        u_1^{(1)}(x)=u_2^{(1)}(x) &\quad \text{ on }\partial\mathcal{S}_h\backslash\partial B_h.
    \end{cases} 
\end{equation}
Then, once again by Theorem \ref{AuxThm}, if $F^{(1)}_1(x_c)\neq F^{(1)}_2(x_c)$, we have that $m^{(1)}(x_c)=0$. This once again contradicts the choice of $\varphi$, so it must be that $F^{(1)}_1(x_c)=F^{(1)}_2(x_c)$ and we have uniquely determined $F^{(1)}$.

Next, we assume that $F^{(1)}_1= F^{(1)}_0$ in $\mathcal{S}_h$ but $F^{(2)}_1\neq F^{(2)}_0$ in $\mathcal{S}_h$, i.e. when $\ell=2$ in Definition \ref{FDef}. In this case, from \eqref{Linear1m}, when $\mathcal{M}_{\omega_1,H,F_1}=\mathcal{M}_{\omega_2,H,F_2}$, in addition to the result on $m^{(1)}$ in the case of $\ell=1$, we have $u^{(1)} := u^{(1)}_0 + (u^{(1)}_1 - u^{(1)}_0)\chi_\omega \in C^{2+\alpha}(\Omega)$ with $u^{(1)}_0 = u^{(1)}_1$ on $\partial\omega$. Then, from \eqref{Linear2}, we can see that $m^{(2)} := m^{(2)}_0 + (m^{(2)}_1 - m^{(2)}_0)\chi_\omega \in C^{2+\alpha}(\Omega)$ with $m^{(2)}_0 = m^{(2)}_1$ on $\partial\omega$ when $\mathcal{M}_{\omega_1,H,F_1}=\mathcal{M}_{\omega_2,H,F_2}$, while $u^{(2)}_j$ satisfies
\begin{equation}
    \begin{cases}
		-D \Delta u^{(2)}_1(x)+\kappa(x)|\nabla u^{(1)}(x)|^2= F^{(1)}(x)m^{(2)}(x)+F^{(2)}_1(x)[m^{(1)}(x)]^2  & \quad\text{ in }  \mathcal{S}_h,\\
        -D \Delta u^{(2)}_0(x)+\kappa(x)|\nabla u^{(1)}(x)|^2= F^{(1)}(x)m^{(2)}(x)+F^{(2)}_0(x)[m^{(1)}(x)]^2  & \quad\text{ in }  B_h\backslash\mathcal{S}_h,\\
        \partial_\nu u_1^{(2)}(x)=\partial_\nu u_0^{(2)}(x) &\quad \text{ on }\partial\mathcal{S}_h\backslash\partial B_h,\\
        u_1^{(2)}(x)=u_0^{(2)}(x) &\quad \text{ on }\partial\mathcal{S}_h\backslash\partial B_h.
    \end{cases}
\end{equation}
Since $u^{(1)}, m^{(1)}, m^{(2)}$ are fixed and smooth, this is simply an elliptic equation. Therefore, we can once again apply the unique continuation principle for elliptic equations \cite{KochTataru2001UCPElliptic} to obtain 
\begin{equation}\label{Linear2m}
    \begin{cases}
		-D \Delta u^{(2)}_1(x)+\kappa(x)|\nabla u^{(1)}(x)|^2 - F^{(1)}(x)m^{(2)}(x) - F^{(2)}_1(x)[m^{(1)}(x)]^2 = 0  & \quad\text{ in }  \mathcal{S}_h,\\
        -D \Delta u^{(2)}_0(x)+\kappa(x)|\nabla u^{(1)}(x)|^2 - F^{(1)}(x)m^{(2)}(x) - F^{(2)}_0(x)[m^{(1)}(x)]^2 =0 & \quad\text{ in }  \mathcal{S}_h,\\
        \partial_\nu u_1^{(2)}(x)=\partial_\nu u_0^{(2)}(x) &\quad \text{ on }\partial\mathcal{S}_h\backslash\partial B_h,\\
        u_1^{(2)}(x)=u_0^{(2)}(x) &\quad \text{ on }\partial\mathcal{S}_h\backslash\partial B_h.
    \end{cases}
\end{equation}
This is simply \eqref{AuxThmEq} with $v_j=u^{(2)}_j$, $P=\kappa|\nabla u^{(1)}|^2 - F^{(1)}m^{(2)}$, $q_j=F^{(2)}_j$ and $h=[m^{(1)}]^2$. Since $F\in\mathcal{D}$ and $m^{(1)}\in C^{2+\alpha}(\mathcal{S}_h)$ by Theorem \ref{ForwardResultStat}, the assumptions of Theorem \ref{AuxThm} are satisfied with $F^{(2)}_1(x_c)\neq F^{(2)}_0(x_c)$, and we have that $h(x_c)=[m^{(1)}(x_c)]^2=0$, i.e. $m^{(1)}(x_c)=0$, for the apex $x_c$ of $\mathcal{S}$. But $m^{(1)}\neq0$ on $\partial\omega_1\supset\partial\mathcal{S}$ by our choice of $\varphi(x)$. Thus, we arrive at a contradiction, and it must be that $\omega_1=\omega_2$ and we have uniquely determined $\omega$.

Moreover, when $\mathcal{M}_{\omega_1,H,F_1}=\mathcal{M}_{\omega_2,H,F_2}$, we have that the solution sets $(u_1,m_1),(u_2,m_2)$ satisfy \eqref{PfF12}.
Repeating the argument above, we have that $m^{(1)}_1 = m^{(1)}_2$, $u^{(1)}_1 = u^{(1)}_2$ and $m^{(2)}_1 = m^{(2)}_2$ in $\mathcal{S}_h$, and $u^{(2)}_j$ satisfies
\begin{equation}
    \begin{cases}
		-D \Delta u^{(2)}_1(x)+\kappa(x)|\nabla u^{(1)}(x)|^2 - F^{(1)}(x)m^{(2)}(x) - F^{(2)}_1(x)[m^{(1)}(x)]^2 = 0  & \quad\text{ in }  \mathcal{S}_h,\\
        -D \Delta u^{(2)}_2(x)+\kappa(x)|\nabla u^{(1)}(x)|^2 - F^{(1)}(x)m^{(2)}(x) - F^{(2)}_2(x)[m^{(1)}(x)]^2 =0 & \quad\text{ in }  \mathcal{S}_h,\\
        \partial_\nu u_1^{(2)}(x)=\partial_\nu u_0^{(2)}(x) &\quad \text{ on }\partial\mathcal{S}_h\backslash\partial B_h,\\
        u_1^{(2)}(x)=u_0^{(2)}(x) &\quad \text{ on }\partial\mathcal{S}_h\backslash\partial B_h.
    \end{cases}
\end{equation}
Then, once again by Theorem \ref{AuxThm}, if $F^{(2)}_1(x_c)\neq F^{(2)}_2(x_c)$, we have that $m^{(1)}(x_c)=0$. This once again contradicts the choice of $\varphi$, so it must be that $F^{(2)}_1(x_c)=F^{(2)}_2(x_c)$ and we have uniquely determined $F^{(2)}$.

The case for higher order $\ell$ follows similarly, by first deducing the equality of $u^{(\ell')}_1=u^{(\ell')}_0$ for $\ell'<\ell$ and $m^{(\ell'')}_1=m^{(\ell'')}_0$ for $\ell''\leq \ell$ on $\partial\omega$, which then gives an elliptic equation in the $\ell$-th order linearisation for $u^{(\ell)}_j$. Then, applying the unique continuation principle for elliptic equations, we obtain \eqref{AuxThmEq} with sufficient regularity of the coefficients. Applying Theorem \ref{AuxThm} then gives $h(x_c)=[m^{(1)}(x_c)]^\ell=0$, which contradicts our choice of $\varphi$. Therefore we have uniquely determined $\omega$. 

Then, the solution sets $(u_1,m_1),(u_2,m_2)$ once again satisfy \eqref{PfF12}. Repeating the argument, we have that $u^{(\ell')}_1=u^{(\ell')}_2$ for $\ell'<\ell$ and $m^{(\ell'')}_1=m^{(\ell'')}_2$ for $\ell''\leq \ell$ in $\mathcal{S}_h$, while $u^{(\ell)}_j$ satisfies an equation of the form \eqref{AuxThmEq}. Applying Theorem \ref{AuxThm} then gives $m^{(1)}(x_c)=0$, which again contradicts our choice of $\varphi$. Thus we have uniquely determined $F^{(\ell)}$. 

\end{proof}

This extends Theorem 2.2 of \cite{LiuMouZhang2022InversePbMeanFieldGames} and Theorem 2.1 of \cite{LiuZhangMFG3} to more general running costs $F$ which are discontinuous across the boundary $\partial\omega$.

		\medskip 
	
	\noindent\textbf{Acknowledgment.} 

	The work was supported by the Hong Kong RGC General Research Funds (projects 11311122, 11300821 and 12301420),  the NSFC/RGC Joint Research Fund (project N\_CityU101/21), and the ANR/RGC Joint Research Grant, A\_CityU203/19.

\bibliographystyle{plain}
\bibliography{ref.bib,refInversePb,refMFG.bib}
\end{document}